\documentclass{amsart}
\usepackage{amssymb,amscd,amsthm,amsmath,amsfonts}

\newcommand{\x}{\times}
\newcommand{\cs}{\mbox{$C^{*}$-algebra}}
\newcommand{\css}{\mbox{$C^{*}$-algebras}}

\newcommand{\C}{\mathbb{C}}

\newcommand{\R}{\mathbb{R}}
\newcommand{\cS}{\mathbb{S}}

\newcommand{\al}{\mbox{$\alpha$}}

\newcommand{\bt}{\mbox{$\beta$}}
\newcommand{\ga}{\mbox{$\gamma$}}
\newcommand{\Ga}{\mbox{$\Gamma$}}
\newcommand{\de}{\mbox{$\delta$}}

\newcommand{\la}{\mbox{$\lambda$}}

\newcommand{\La}{\mbox{$\Lambda$}}

\newcommand{\Om}{\mbox{$\Omega$}}

\newcommand{\mfA}{\mathfrak{A}}

\newcommand{\bgc}{\begin{center}}

\newcommand{\edc}{\end{center}}
\newcommand{\be}{\begin{enumerate}}
\newcommand{\ee}{\end{enumerate}}
\newcommand{\beq}{\begin{equation}}
\newcommand{\eeq}{\end{equation}}
\newcommand{\beqn}{\begin{eqnarray}}
\newcommand{\eeqn}{\end{eqnarray}}
\newcommand{\beqns}{\begin{eqnarray*}}
\newcommand{\eeqns}{\end{eqnarray*}}
\newcommand{\bq}{\begin{quote}}
\newcommand{\eq}{\end{quote}}

\newcommand{\bi}{\begin{itemize}}
\newcommand{\ei}{\end{itemize}}
\newcommand{\bd}{\begin{description}}
\newcommand{\ed}{\end{description}}

\theoremstyle{plain}
\newtheorem{theorem}{Theorem}

\newtheorem{definition}{Definition}
\newtheorem{proposition}{Proposition}
\newtheorem{corollary}{Corollary}

\numberwithin{equation}{section}

\begin{document}
\title{Continuous family groupoids}
\author{Alan L. T. Paterson}
\address{Department of Mathematics\\
University of Mississippi\\
University, MS 38677}
\email{mmap@olemiss.edu}
\keywords{Groupoids, continuous families, continuous left Haar systems,
$G$-spaces, index theorems}
\subjclass{Primary: 22A22; 58H05; Secondary 58G10}
\date{}
\begin{abstract}
In this paper, we define and investigate the properties of continuous family
groupoids. This class of groupoids is necessary for investigating the groupoid
index theory arising from the equivariant Atiyah-Singer index theorem for
families, and is also required in noncommutative geometry. The class includes
that of Lie groupoids, and the paper shows that, like Lie groupoids, continuous
family groupoids always admit (an essentially unique) continuous left Haar
system of smooth measures. We also show that the action of a continuous family
groupoid $G$ on a continuous family $G$-space fibered over another continuous
family $G$-space $Y$ can always be regarded as an action of the continuous
family groupoid $G*Y$ on an ordinary $G*Y$-space.
\end{abstract}

\maketitle

\section{Introduction}
In the context of noncommutative geometry, it is becoming increasingly clear
that the classical Atiyah-Singer index theory, in which a compact group acts
equivariantly on a compact manifold, requires to be extended to the context
of a Lie groupoid acting properly on a manifold.   In this
connection, Connes (\cite[p.151]{Connesbook}) says:
\bq
One of the interests of the general formulation [of the Baum-Connes
conjecture] is to put many particular results in a common framework.  Thus,
for instance, the following three theorems:
\bi
\item[1)] The Atiyah-Singer index theorem for covering spaces
\item[2)] The index theorem for measured foliations
\item[3)] The index theorem for homogeneous spaces
\ei
are all special cases of the same index theorem for $G$-invariant elliptic
operators $D$ on proper $G$-manifolds, where $G$ is a smooth [i.e. Lie]
groupoid with a transverse measure $\La$.
\eq
For example, in the foliation case, the smooth (Lie) groupoid $G$ is the
holonomy groupoid of a compact foliated manifold. Like all Lie groupoids, it
acts properly on itself as a $G$-space, each $g\in G$ acting as a
diffeomorphism from $G^{s(g)}$ onto $G^{r(g)}$. (Here, $s,r$ are respectively
the source and range maps on the groupoid. See \S 2 for more details.)
So the ``proper $G$-manifold'' in that case is $G$ itself. Each $G$-invariant
family of elliptic pseudodifferential operators along the leaves defines an
element of $K^{0}(F^{*})$, where $F$ is the vector bundle of vectors tangential
along the $G^{u}$'s, and has its index in $K_{0}(C^{*}(G))$.

We note that smoothness is really only being used along the ``leaves''
$G^{u}$ (and $G_{u}$), and this leads naturally to considering the case
where only this kind of smoothness is assumed for the groupoid, and we have
only continuity ``transversely''. (This notion, as well as the
$C^{0,\infty}$-notation of the paper, appears in \cite{Cointeg}.)
In fact, such a class of groupoids is required for the equivariant version of
the Atiyah-Singer index theorem for continuous families (\cite{AS4}) alluded
to at the end of that paper (\cite[p.135]{AS4}). (The groupoid interpretation
of this theorem is given in \S 4 of the present paper.) This equivariant index
theorem cannot be formulated in terms of Lie groupoids since the unit space of
the groupoid in
that context is the base space $Y$ of the continuous family on which a
compact Lie group acts, and $Y$ is not assumed to be a manifold but only a
compact Hausdorff space. To deal with this, then, we need to consider a class
of locally compact groupoids more general than that of Lie groupoids, in which,
as above, we only have smoothness along the leaves. The groupoids that we need
for this are the {\em continuous family groupoids} in the title of the paper.

Continuous families (in the sense of Atiyah and Singer) play an important role
in index theory. For example, the Bott periodicity theorem, which asserts (for
compact $Y$) that $K^{0}(Y)=K^{0}(Y\x \cS^{2})$, involves the continuous
family $Y\x \cS^{2}$ over $Y$.  So if we restricted only to smooth families,
we would be effectively restricting K-theory to smooth manifolds, whereas
K-theory is a functor on the topological category.

In the non-equivariant families theorem, a continuous
family is defined as follows. We are given a compact Hausdorff space $Y$ and a
fiber bundle $X$ over $Y$ with fiber $Z$, where $Z$ is a smooth compact
manifold and the structure group of $X$ is $\text{Diff}(Z)$.\footnote{In the
notation of \cite{AS4}, the $X$ and $Z$ are interchanged.} (Here
$\text{Diff}(Z)$ is a topological group under the topology of uniform
convergence for each derivative.) The space $X$ is thus a continous family of
spaces diffeomorphic to $Z$. The index of an elliptic family on $X$ is shown to
lie in $K^{0}(Y)$.

For continuous family groupoids, we need to extend this to the case where the
fibers are not diffeomorphic to a fixed space. The general notion of a
continuous family of manifolds, required for the paper, is more accurately
described as being a {\em locally} continuous family $X$ of smooth manifolds
$X^{y}$ with $y\in Y$, but for the sake of brevity, we will omit the adjective
``locally''. (The smooth version of this for almost differentiable groupoids is
given in \cite{NWX}.) The role of $\text{Diff}(Z)$ is taken over by a certain
pseudogroup of maps. We then (as in \cite{AS4}) describe what we mean by a
vector bundle over $X$ which is smooth along the fibers. (This is not actually
used later in the paper but is included since it is required for groupoid index
theory.)

The definition of {\em continuous family groupoids} is given in \S 3. {\em
Lie groupoids} are, of course, continuous family groupoids. But there
are many naturally occurring examples of continuous family groupoids
that are not Lie groupoids (including, in particular, the groupoid associated
with the equivariant index theorem for families referred to earlier). As for
Lie groupoids, there is an essentially unique, continuously varying, left
Haar system of smooth measures
for any continuous family groupoid $G$, so that there is a
canonical $C^{*}(G)$. (The $K_{0}$-group of this \cs\ is the recipient for the
index of elliptic families on $G$-manifolds.)

The last section \S 4 of the paper discusses $G$-spaces.  Its main result
is that for the category of continuous
family groupoids, by changing the groupoid, we need
only consider ``ordinary'' $G$-spaces rather than fibered $G$-spaces. Indeed,
suppose that we are given $G$-spaces $X,Y$ with $X$ fibered equivariantly over
$Y$. Then we can form the ``transformation groupoid'' $G*Y$ associated with the
action of $G$ on $Y$. This groupoid is shown to be itself a continuous family
groupoid, and the action of $G$ on $X$ fibered over $Y$ is equivalent to the
canonical action of $G*Y$ on $X$ with $X$ treated as an ordinary $G*Y$-space.

In the case of the equivariant Atiyah-Singer families index theorem, in which
a compact Lie group $H$ acts equivariantly on a compact manifold $X$ over $Y$,
the preceding shows that this is equivalent to the transformation group
groupoid $H*Y$ acting on $X$, and even in that classical context, we leave the
group category for the continuous family groupoid category. (The index of a
$H$-invariant elliptic family of pseudodifferential operators on $X$ can be
shown to lie in $K_{0}(C^{*}(H*Y))$.)

\section{Continuous families of manifolds}
Let $Y$ be a topological space and $k\geq 1$. Let $A_{1},A_{2}$ be open subsets
of $Y\x \R^{k}$ such that $q_{1}(A_{1})\subset q_{1}(A_{2})$ where $q_{1}$ is
the canonical projection map from $Y\x \R^{k}$ onto $Y$. For $y\in
q_{1}(A_{i})$ ($i=1,2$), let $A_{i}^{y}=\{x\in \R^{k}: (y,x)\in A_{i}\}$. Let
$f:A_{1}\to A_{2}$ be a continuous
map which preserves fibers, i.e. for each
$y\in q_{1}(A_{1})$, we have $f(\{y\}\x A_{1}^{y})\subset\{y\}\x A_{2}^{y}$. For
such an $f$, define $f^{y}:A_{1}^{y}\to A_{2}^{y}$ by: $f(y,x)=(y,f^{y}(x))$.

Then (cf. \cite[p. 110]{Cointeg}) the function $f$ is said to be a {\em
$C^{0,\infty}$-function ($f\in C^{0,\infty}(A_{1},A_{2})$)} if,
whenever $U_{1}, U_{2}$ are open
subsets of $Y$ and $V_{1}, V_{2}$ are open subsets of $\R^{k}$ such that
$U_{i}\x V_{i}\subset A_{i}$ for each $i$
and $f(U_{1}\x V_{1})\subset U_{2}\x V_{2}$, then
the  map $y\to (f^{y})_{\mid V_{1}}$ is a continuous  map
from $U_{1}$ into $C^{\infty}(V_{1}, V_{2})$. Here, the topology on
$C^{\infty}(V_{1}, V_{2})$ is that of uniform convergence on compacta for all
derivatives, i.e. $h_{n}\to h$ in $C^{\infty}(V_{1}, V_{2})$ if and only if,
for every compact subset $K$ of $V_{1}$ and multi-index
$\al=(\al_{1}, \al_{2},\ldots ,\al_{k})$, we have
\begin{equation}
\left|\partial^{\al}h_{n}(x)- \partial^{\al}h(x)\right|_{K}\to 0
\label{eq:par}
\end{equation}
as $n\to \infty$. Here, the differentiation is with respect to the $x_{i}$'s
where $x=(x_{1},\ldots ,x_{k})$, and
$\left| g\right|_{K}=\sup_{x\in K}\left| g(x)\right|$ for any complex-valued
function $g$ bounded on $K$.

In the above definition, since $f(U_{1}\x V_{1})\subset U_{2}\x
V_{2}$ and $f$ is fiber preserving, we have $U_{1}\subset U_{2}$.
We can clearly take $U_{1}=U_{2}$ above.  In addition, by
the continuity of $f$, every element of $A_{1}$ belongs to {\em some} open
$U_{1}\x V_{1}\subset A_{1}$ for which there exists an open $U_{2}\x
V_{2}\subset A_{2}$ with $f(U_{1}\x V_{1})\subset U_{2}\x V_{2}$.

The set $\text{Diff}^{0}(A_{1},A_{2})$ is defined to be the set of functions
$f\in C^{0,\infty}(A_{1},A_{2})$ for which $f^{-1}$ exists and belongs to
$C^{0,\infty}(A_{2},A_{1})$.  (If $\text{Diff}^{0}(A_{1},A_{2})$
is non-empty, then of course $q_{1}(A_{1})=q_{1}(A_{2})$.) Every element of
$\text{Diff}^{0}(A_{1},A_{2})$ is trivially a homeomorphism.

In practice, an alternative formulation of $C^{0,\infty}(A_{1},A_{2})$
proves useful.  Let us say that a fiber preserving function $f:A_{1}\to
A_{2}$ is $C^{\infty}${\em -continuous} if given $a\in A_{1}$ and an
open subset of $A_{2}$ of the form $U_{2}\x V_{2}$ ($U_{2}\subset Y,
V_{2}\subset \R^{k}$) which contains $f(a)$, then there exists an open subset
$U_{1}\x V_{1}$ of $A_{1}$ such that $a\in U_{1}\x V_{1}$, $f(U_{1}\x
V_{1})\subset U_{2}\x V_{2}$ and the map $y\to (f^{y})_{\mid V_{1}}$ takes
$U_{1}$ into $C^{\infty}(V_{1}, V_{2})$.

\begin{proposition}       \label{prop:cinf}
A function $f:A_{1}\to A_{2}$ belongs to $C^{0,\infty}(A_{1},A_{2})$ if and
only if it is $C^{\infty}$-continuous.
\end{proposition}
\begin{proof}
If $f\in C^{0,\infty}(A_{1},A_{2})$, then trivially $f$ is
$C^{\infty}$-continuous.  Conversely, suppose that $f$ is
$C^{\infty}$-continuous and let $U_{i}\x V_{i}$ be open subsets of $A_{i}$
such that $f(U_{1}\x V_{1})\subset U_{2}\x V_{2}$.
Let $K$ be a compact subset of $V_{1}$ and $y_{0}\in U_{1}$.
Then for each $v\in K$, there exists, by the $C^{\infty}$-continuity
of $f$, an open subset $U_{i}^{v}\x V_{i}^{v}$ of $A_{i}$ with
$(y_{0},v)\in U_{1}^{v}\x V_{1}^{v}$ such that the map
$y\to (f^{y})_{\mid V_{1}^{v}}$ is continuous from $U_{1}^{v}$ into
$C^{\infty}(V_{1}^{v},V_{2}^{v})$.  Cover $K$ by a finite number of sets
$V_{1}^{v_{1}},\ldots ,V_{1}^{v_{n}}$, and let
$U'=\cap_{j=1}^{n} U_{1}^{v_{j}}$.  There exist compact subsets $K^{j}$
$(1\leq j\leq n)$ such that $K^{j}\subset V_{1}^{v_{j}}$ and
$K=\cup_{j=1}^{n} K^{j}$. Then each
of the maps $y\to (f^{y})_{\mid V_{1}^{v_{j}}}$ is continuous, and it
follows (using the $K^{j}$'s) that
$\left\|\partial^{\al}f^{y} - \partial^{\al}f^{y_{0}}\right\|_{K}\to 0$ if
$y\to y_{0}$ in $U'$.  So $f\in C^{0,\infty}(A_{1},A_{2})$.
\end{proof}
\begin{corollary}  \hspace{.1in}        \label{cor:cinf}
Let $A_{1},A_{2}, A_{3}$ be open subsets
of $Y\x \R^{k}$ such that $q_{1}(A_{i})\subset q_{1}(A_{i+1})$ $(1\leq i\leq
2)$.  Let $f\in C^{0,\infty}(A_{1},A_{2})$ and
$g\in C^{0,\infty}(A_{2},A_{3})$.  Then
$g\circ f\in C^{0,\infty}(A_{1},A_{3})$.  Further, if
$f\in \text{Diff}^{0}(A_{1},A_{2})$, then
$f^{-1}\in \text{Diff}^{0}(A_{1},A_{2})$.
\end{corollary}
\begin{proof}
For the first assertion of the corollary, one just has to prove that $g\circ f$
is $C^{\infty}$-continuous. To this end, note that $(g\circ f)^{y}=g^{y}\circ
f^{y}$. One then follows the elementary proof that
the composition of two continuous functions is continuous, and uses induction
and the chain rule to deal with the partial derivatives in (\ref{eq:par}). The
second assertion is obvious.
\end{proof}

We now recall the definition of a {\em pseudogroup}.  Various definitions
have been given of this in the literature: the version that we will use is
that given in \cite[p.1]{KobNom}.  A {\em pseudogroup} $S$ on
a topological space $X$ is an inverse semigroup of homeomorphisms
$f:A\to B$, where $A, B$ are open subsets of $X$ (depending on $f$) such that:
\bi
\item[(i)] if $f:A\to B$ is a homeomorphism, where $A=\cup_{i\in I} A_{i}$ and
$f_{\mid A_{i}}\in S$, then $f\in S$;
\item[(ii)] if $f:A\to B$ belongs to $S$ and $A'$ is an open subset
of $A$, then $f_{\mid A'}\in S$.
\item[(iii)] If $A$ is open in $X$, then the identity map
$\text{id}:A\to A$ belongs to $S$.
\ei

Let $\Gamma(Y\x \R^{k})$ be the union of all of the sets
$\text{Diff}^{0}(A_{1},A_{2})$ (with $A_{1}, A_{2}$ ranging over the open
subsets of $Y\x \R^{k}$).

\begin{proposition}  \hspace{.1in}        \label{prop:gamma}
The set $\Gamma(Y\x \R^{k})$ is a pseudogroup on $Y\x \R^{k}$.
\end{proposition}
\begin{proof}
Let $S=\Gamma(Y\x \R^{k})$.  By Corollary~\ref{cor:cinf}, if $f\in
\text{Diff}^{0}(A_{1},A_{2})$ and  $g\in \text{Diff}^{0}(A_{3},A_{4})$,
then $f^{-1}\in S$ and both $f\circ g$,
$(f\circ g)^{-1}=g^{-1}\circ f^{-1}$ are $\text{Diff}^{0}$ maps.
From Corollary~\ref{cor:cinf}, if $f\in \text{Diff}^{0}(A_{1},A_{2})$ and
$g\in \text{Diff}^{0}(A_{3},A_{4})$ then
$g\circ f\in \text{Diff}^{0}(A_{5},A_{6})$, where
$A_{5}=f^{-1}(A_{2}\cap A_{3})$ and $A_{6}=g(A_{2}\cap A_{3})$. So
$f^{-1}, f\circ g\in S$, and $S$ is an inverse semigroup.  Conditions (i)
and (ii) above follow by $C^{\infty}$-continuity (cf. the proof of
Proposition~\ref{prop:cinf}), while (iii) is trivial.
\end{proof}

Let $X$, $Y$ be locally compact Hausdorff spaces and $p:X\to Y$ be a
continuous open surjection. We say that {\em $(X,p)$ is fibered over $Y$}
with fibers $X^{y}=p^{-1}(\{y\})$ ($y\in Y$), and call $(X,p)$ a {\em fiber
space} (over $Y$). We now define what is meant by a continuous family of
manifolds over $Y$.

\begin{definition}             \label{def:cf}
Let $(X,p)$ be fibered over $Y$.  Then the pair $(X,p)$
is defined to be a {\em continuous family of manifolds $X^{y}$ over $Y$} or
simply {\em a continuous family over $Y$} if
there exists a set of pairs $\{(U_{\al},\phi_{\al}): \al\in A\}$, where each
$U_{\al}$ is an open subset of $X$ and $\cup_{\al\in A}U_{\al}=X$,
compatible with the pseudogroup $\Gamma(Y\x \R^{k})$ in the following sense:
\bi
\item[(i)] for each $\al$, the map
$\phi_{\al}$ is a homeomorphism from $U_{\al}$ onto an open subset of
$Y\x \R^{k}$ for which $q_{1}\circ \phi_{\al}=p_{\mid U_{\al}}$;
\item[(ii)] for each $\al, \bt$, the mapping
$\phi_{\bt}\circ\phi_{\al}^{-1}\in
\text{Diff}^{0}(\phi_{\al}(U_{\al}\cap U_{\bt}),
\phi_{\bt}(U_{\al}\cap U_{\bt})).$
\ei
\end{definition}

The family $\mfA=\{U_{\al}:\al\in A\}$ will be called an {\em atlas} for the
continuous family $(X,p)$, and the $U_{\al}$'s, or more precisely, the pairs
$(U_{\al},\phi_{\al})$, will be called {\em charts}.
Of course, in the above definition, we can and will take the atlas $\mfA$ to be
maximal. Then $\mfA$ is a basis for the topology of $X$.

If $(U,\phi)\in \mfA$ and $x\in U$, then there exists a $V\subset U$
with $x\in V$ such that $(V,\phi\mid_{V})\in \mfA$ and $\phi(V)=p(V)\x W$ for
some open subset $W$ of $\R^{k}$.805We shall write $V\sim p(V)\x W$.  For many
purposes, we need only consider charts of this special form $V$.

The simplest example of a continuous family over $Y$ is one of the form
$X=Y\x M$ where $M$ is a manifold. Such a family is called {\em trivial}.
From the preceding paragraph, every continuous family is locally trivial.

A continuous family in the sense of Atiyah and Singer is a continuous family in
our sense. To see this, recall (\S 1)
that in that case, $(X,p)$ is a fiber bundle
over $Y$ with a manifold $Z$ as fiber and with structure group
$\text{Diff}(Z)$. Then there is a basis $\{J_{\de}\}$ for the topology of $Y$
and fiber preserving homeomorphisms $c_{\de}:p^{-1}(J_{\de})\to J_{\de}\x Z$
such that the resultant cocycles are continuous maps into $\text{Diff}(Z)$.
Let $(L_{\ga},\chi_{\ga})$ be a chart for $Z$. We obtain charts for the
continuous family in the sense of Definition~\ref{def:cf} by taking sets of the
form $((1\otimes \chi_{\ga})\circ c_{\de})^{-1}(A)$ where $A$ is an open subset
of $J_{\de}\x \R^{k}$ ($k=\dim Z$)..

{\em Smooth} families are, of course, continuous families. These arise in the
theory of Lie groupoids (\S 3) (\cite{Mackenzie,NWX,Paterson,Pra66,Pra67}). (In
particular, for any Lie groupoid $G$ with range and source maps $r,s$ and unit
space $G^{0}$, both $(G,r), (G,s)$ are smooth families over $G^{0}$.) For a
smooth family, we require that both $X,Y$ be manifolds and that $p:X\to Y$ be
a (surjective) submersion. Then locally, $p$ can be taken to be a smooth
projection map $(x,y)\to x$, and thus
defines a foliated manifold structure on $X$
whose leafs are the $X^{y}$'s (\cite[p.23-24]{Camacho}). Condition (ii) of
Definition~\ref{def:cf} is satisfied since the maps
$\phi_{\bt}\circ\phi_{\al}^{-1}$ are diffeomorphisms.

Let $(X,p)$ be a continuous family of manifolds with atlas
$\mfA=\{U_{\al}: \al\in A\}$. Then (as is to be expected)
every $X^{y}$ is a (smooth) manifold.  Indeed, for fixed $y$, let
\[        \mfA^{y}=\{U_{\al}\cap X^{y}: \al\in A\}.           \]
Then $\mfA^{y}$ gives the relative topology on $X^{y}$ (as a closed subset
of $X$).  Since the restriction of $\phi_{\bt}\circ\phi_{\al}^{-1}$ to
$\phi_{\al}(U_{\al}\cap U_{\bt}\cap X^{y})$ is a diffeomeomorphism onto
$\phi_{\bt}(U_{\al}\cap U_{\bt}\cap X^{y})$, we obtain that $X^{y}$ is a
manifold.

We now describe some operations that produce new continuous families from
given ones. We note without giving details that similar (easier) constructions
can be given for fiber spaces.

A pull-back of a continuous family is also a continuous family.
Specifically, let $(X,p)$ be a continuous family over $Y$, $Z$ be a locally
compact Hausdorff space and $t:Z\to Y$ be a continuous map.  The pull-back
continuous family $(t^{-1}X,p')$ over $Z$ is given by the subset
\[            t^{-1}X=\{(z,x)\in A\x X: t(z)=p(x)\}               \]
of $Z\x X$ and the map $p'$ where $p'((z,x))=z$. We have the commuting diagram:
\begin{equation} \label{CD:pullback}
\begin{CD}
	t^{-1}X     @>t'>>    X  \\
      @Vp'VV               @VVpV \\
	 Z         @>t>>      Y
\end{CD}
\end{equation}
where $t'((z,x))=x$.  We obtain charts for $(t^{-1}X,Z)$ as follows.
If $(U,\phi)$ is a chart for $X$ and
$\phi(x)=(p(x),h(x))$, then $(t^{-1}U,\phi')$ is a chart for $t^{-1}X$, where
\[       \phi'((z,x))= (z,h(x)).               \]

Let $(X_{1},p_{1}), (X_{2},p_{2})$ be continuous families over the same space
$Y$ and $f:X_{1}\to X_{2}$ be a continuous fiber preserving map, i.e.
$p_{2}\circ f=p_{1}$. We say that $f\in C^{0,\infty}(X_{1},X_{2})$ if $f$ is
locally $C^{0,\infty}$ in the earlier sense. That is, whenever $(U,\phi)$,
$(V,\psi)$ are charts for $X_{1}$ and $X_{2}$ respectively such that
$p_{1}(U)=p_{2}(V)$ and $f(U)\subset V$, then $\psi\circ f\circ
\phi^{-1}\in C^{0,\infty}(\phi(U),\psi(V))$.

We now discuss what is meant by a {\em morphism} of continuous families.  We
deal first with a special case (to which, as we shall see, the general case can
be reduced).  In the special case, a {\em morphism} from $(X_{1},p_{1})$
into $(X_{2},p_{2})$  over $Y$
as in the preceding paragraph is just a function $f\in
C^{0,\infty}(X_{1},X_{2})$. We represent such a morphism by
the commutative diagram:
\begin{equation} \label{CD:morphismspecial}
\begin{CD}
	X_{1}     @>f>>    X_{2}  \\
      @Vp_{1}VV @VVp_{2}V \\
	     Y @>=>> Y
\end{CD}
\end{equation}

We now define what is meant by a morphism of continuous families in the general
case. Let $(X_{1},p_{1}), (X_{2},p_{2})$ be continuous families over locally
compact Hausdorff spaces $Y_{1}, Y_{2}$. Let $q:Y_{1}\to Y_{2}$ be a continuous
map and $f:X_{1}\to X_{2}$ be a continuous fiber preserving map with respect
to $q$ in the sense that $p_{2}\circ f=q\circ p_{1}$. (In the special case of a
morphism above, $Y_{1}=Y_{2}=Y$ and $q$ is the identity map.) For $y\in Y_{1}$,
let $f^{y}:X_{1}^{y}\to X_{2}^{q(y)}$ be the restriction of $f$ to $X_{1}^{y}$.
Then $f$ is called a {\em morphism} (with respect to $q$) if for each $y\in
Y_{1}$, the function $f^{y}\in C^{\infty}(X_{1}^{y},X_{2}^{q(y)})$
and the map $y\to
f^{y}$ is locally a $C^{0,\infty}$-function. More precisely, given $x\in
X_{1}$ and a chart $U_{2}\sim p_{2}(U_{2})\x W_{2}$ in $X_{2}$ containing
$f(x)$, then there exists a chart $U_{1}\sim
p_{1}(U_{1})\x W_{1}$ in $X_{1}$, and such that
$x\in U_{1}$, $f(U_{1})\subset U_{2}$
and (in local coordinates) the map $y\to f^{y}$ is continuous from
$p_{1}(U_{1})$ into $C^{\infty}(W_{1},W_{2})$ where
\[          f((y,w_{1}))=(q(y),f^{y}(w_{1})).                 \]
We think of the continuous family $X_{1}$ over $Y_{1}$ as being taken over into
the continuous family $X_{2}$ over $Y_{2}$ by the maps $f,q$.

We represent a morphism $f$ by the following commutative diagram:
\begin{equation} \label{CD:morphism}
\begin{CD}
	X_{1}     @>f>>    X_{2}  \\
      @Vp_{1}VV @VVp_{2}V \\
      Y_{1} @>q>> Y_{2}
\end{CD}
\end{equation}
The set of morphisms $f$ from $X_{1}$ into $X_{2}$ is denoted by
$C^{0,\infty}(X_{1},X_{2})$.

It is easy to
prove that if $(X_{3},p_{3})$ is a continuous family over $Y_{3}$,
$q':Y_{2}\to Y_{3}$ is
continuous and $g:X_{2}\to X_{3}$ is a morphism with respect to $q'$, then
$g\circ f:X_{1}\to X_{3}$ is a morphism with respect to $q'\circ q$.  The
commuting diagram for the morphism $g\circ f$ can be represented as the
``product'' of the following commuting diagram:
\[
\begin{CD}
     X_{1}  @>>f>     X_{2}   @>>g>     X_{3} \\
	@Vp_{1} VV    @Vp_{2} VV         @Vp_{3}VV \\
	Y_{1}  @>q>>     Y_{2}  @>q'>>    Y_{3}
\end{CD}
\]
It is easy to show that the class of continuous families is a category with
morphisms $f$ as above.

A morphism $f:X_{1}\to X_{2}$ in the general case can be reduced simply
to the special case by using a pull-back continuous family. More precisely, the
pull-back continuous family $q^{-1}X_{2}$ is a continuous family over $Y_{1}$,
and we have a map $f':X_{1}\to q^{-1}X_{2}$ given by: $f'(x)=(p_{1}(x),f(x))$.
It is left to the reader to check that a continuous, fiber-preserving map $f$
is a morphism if and only if $f'$ is a (special case) morphism.

In the case where $X=X_{1}, Y=Y_{1},
X_{2}=\C$ and $Y_{2}$ is a singleton, then we write
$C^{0,\infty}(X_{1},\C)=C^{0,\infty}(X_{1})$. A simple
argument (using the charts for $X$) shows that $C^{0,\infty}(X)$ is a
$^{*}$-subalgebra of $C(X)$, the algebra of continuous complex-valued
functions on $X$. It is left to the reader to show that $X$ admits
$C^{0,\infty}-$ partitions of unity.

Let $C_{c}^{0,\infty}(X)$ be the subalgebra of functions with compact support
in $C^{0,\infty}(X)$. If $X$ is a product $Y\x V$ where $V$
is an open subset of $\R^{k}$ and $(y_{0},k_{0})\in X$, then there is a
function $f\in C_{c}^{0, \infty}(X)$ with $f((y_{0},k_{0}))=1$. Indeed, we can
take $f=g\otimes h$ where $g\in C_{c}(Y), h\in C_{c}^{\infty}(V)$ and
$g(y_{0})=1=h(k_{0})$. By considering charts, it follows that for
a general continuous family $X$,
the space $C_{c}^{0,\infty}(X)$ separates the points of $X$.

Now let $(X_{1},p_{1}), (X_{2},p_{2})$ be continuous families
over $Y$. Let $(X_{1}*X_{2},p)$ be the {\em fibered product}
of $(X_{1},p_{1})$ and $(X_{2},p_{2})$: so
\[ X_{1}*X_{2}=\{(x_{1},x_{2})\in X_{1}\x X_{2}: p_{1}(x_{1})
=p_{2}(x_{2})\} \]
and $p(x_{1},x_{2})=p_{1}(x_{1})=p_{2}(x_{2})$. We sometimes write
$p=p_{1}*p_{2}$. Then with the
relative topology, $(X_{1}*X_{2},p)$ is
a continuous family of manifolds over $Y$, with
each fiber $X^{y}={X_{1}}^{y}\x {X_{2}}^{y}$ having the product manifold
structure. Indeed, it is left to the reader to check that if $(x_{1},x_{2})\in
X_{1}*X_{2}$ and $(U_{i},\phi_{i})$ are charts for $X_{i}$ $(i=1,2)$, $x_{i}\in
U_{i}, p_{1}(U_{1})=p_{2}(U_{2})$, then in an obvious notation,
$(U_{1}*U_{2},\phi_{1}*\phi_{2})$ is a chart for $X$ and these charts determine
an atlas for $X_{1}*X_{2}$. Further, $p$ is a continuous surjection which is
open since $p(U_{1}*U_{2})=p_{1}(U_{1})$ is open in $Y$.

For clarity, we will sometimes write $(X_{1},p_{1})*(X_{2},p_{2})$ in place
of $X_{1}*X_{2}$. Note that $Y$ is trivially a continuous family
over itself and that $Y*X_{1}\cong X_{1}$.  Note also that if $X_{1}', X_{2}'$
are continuous families over $Y$ and if $a_{i}:X_{i}'\to X_{i}$ are
morphisms, then the natural map $a_{1}*a_{2}: X_{1}'*X_{2}'\to X_{1}*X_{2}$
is a morphism of continuous families.  This map is an isomorphism if both
$a_{1}, a_{2}$ are.  These are proved by reducing to the case of charts.

We will need a slight generalization of $X_{1}*X_{2}$ above later.  In this
situation, we have, for $i=1,2$, fiber spaces $(Y_{i},v_{i})$ over some
$Z$, and $(X_{i},p_{i})$ continuous families over $Y_{i}$.  So $Y_{1}*Y_{2}$
is a fiber space over $Z$. We can then form
the continuous family $(X_{1}*X_{2},p_{1}*p_{2})$ over $Y_{1}*Y_{2}$, where
\[ X_{1}*X_{2}=\{(x_{1},x_{2}): x_{i}\in X_{i},
p_{1}*p_{2}(x_{1},x_{2})\in Y_{1}*Y_{2}\}     \]
and $p_{1}*p_{2}(x_{1},x_{2})=(p_{1}(x_{1}),p_{2}(x_{2})$. Note that when
$Y_{1}=Y_{2}=Y=Z$ with each $v_{i}$ the identity map, then the two definitions
of $(X_{1}*X_{2},p_{1}*p_{2})$ coincide. It is left to the reader to show
that $(X_{1}*X_{2},p_{1}*p_{2})$ is a continuous family over $Y_{1}*Y_{2}$.

Now let $(X_{1},p_{1})$ and $(X_{2},p_{2})$ be continuous families over $Y$.
There are two other ways in which $X_{1}*X_{2}$ can be naturally regarded as
a continuous family.  These give pull-back families
$(X_{1}*X_{2},t_{1})$ over $X_{1}$ and
$(X_{1}*X_{2},t_{2})$ over $X_{2}$.  Here the $t_{i}$'s are the natural
projection maps:
\[ t_{1}(x_{1},x_{2})=x_{1}, t_{2}(x_{1},x_{2})=x_{2}.\]
Firstly, the pull-back of the continuous family $(X_{2},p_{2})$ by
$p_{1}:X_{1}\to Y$ gives the continuous family
$(p_{1}^{-1}(X_{2}),t_{1})$, where
\[ p_{1}^{-1}(X_{2})=\{(x_{1},x_{2}): p_{1}(x_{1})=p_{2}(x_{2})\}=
X_{1}*X_{2}.\]
In the same way, $(p_{2}^{-1}(X_{1}),t_{2})=
(X_{2}*X_{1},t_{1})$ over $X_{2}$.  Interchanging
first and second components gives the continuous family $(X_{1}*X_{2},t_{2})$.

Note that if $(X_{2},p_{2})$ is assumed to be only a fiber space over $Y$
(with $(X_{1},p_{1})$ still being assumed to be a continuous family), then
$(X_{1}*X_{2},t_{2})$ is still a continuous family, the $X_{2}$ just playing
the continuous role of a parameter space.
The following proposition will be used in \S 4.

\begin{proposition}           \label{prop:morXYZ}
Let $(X_{1},q_{1})$ and $(X_{2},q_{2})$ be continuous families over $Z$ and
$f:X_{1}\to X_{2}$ be a morphism.  Let $(Y,q)$ be a continuous family over $Z$.
Then there is a canonical morphism $f*1$ from $(X_{1}*Y,t_{2})$ into
$(X_{2}*Y,t_{2})$:
\begin{equation} \label{CD:morXYZ}
\begin{CD}
	X_{1}*Y     @>f*1>>    X_{2}*Y  \\
      @Vt_{2}VV             @VVt_{2}V \\
	 Y           @>=>>       Y
\end{CD}
\end{equation}
\end{proposition}
\begin{proof}
To obtain (\ref{CD:morXYZ}), we ``$*$'' the morphism diagram for $f$ by $Y$
to get:
\begin{equation} \label{CD:f*1}
\begin{CD}
	X_{1}*Y     @>f*1>>    X_{2}*Y  \\
      @Vq_{1}*1VV            @VVq_{2}*1V \\
      Z*Y        @>=>>       Z*Y
\end{CD}
\end{equation}
It is easy to check that (\ref{CD:f*1}) is commutative, and it is left to the
reader to check, using charts, that $f*1$ is a morphism. (\ref{CD:morXYZ})
now follows by identifying $Z*Y$ with $Y$ and $q_{i}*1$ with $t_{2}$.
\end{proof}

Now let $(X,p)$ be a continuous family of manifolds over $Y$. We will
define what is meant by a {\em smooth vector bundle over $X$} of dimension $m$.
This generalizes the corresponding notion in \cite{AS4}.

We first define the pseudogroup in this situation corresponding to $\Ga(Y\x
\R^{k})$ earlier.  Let $A_{1}, A_{2}$ be open subsets of $Y\x \R^{k}$ as
earlier.  Let $m\geq 0$ and $B_{i}=A_{i}\x \R^{m}$.  Note that
each $B_{i}\subset Y\x \R^{k+m}$.  Let $f:A_{1}\to A_{2}$ and
$g:B_{1}\to \R^{m}$, and define $F:B_{1}\to B_{2}$ by:
\[         F(y,w,\xi)=(y,f(y,w), g(y,w,\xi))           \]
where $(y,w)\in A_{1}, \xi\in \R^{m}$.
We say that $F\in \text{Diff}^{0}_{\ell}(B_{1},B_{2})$ if $F\in
\text{Diff}^{0}(B_{1},B_{2})$ ($B_{1}, B_{2}$ regarded as subsets of
$Y\x\R^{m}$) and for fixed $(y,w)$, the map
$\xi\to g(y,w,\xi)$ is a vector space isomorphism of $\R^{m}$.  Let
$\Ga(Y\x \R^{k};\R^{m})$ be the union of all of the sets
$\text{Diff}^{0}_{\ell}(B_{1},B_{2})$.  It is easy to check that
$\text{Diff}^{0}_{\ell}(B_{1},B_{2})$ is a subpseudogroup of
$\Ga(Y\x \R^{k+m})$.

\begin{definition}   \label{def:smv}
An $m$-dimensional
vector bundle $(E,\pi)$ over $X$ is said to be a {\em smooth vector bundle}
if there exists an atlas of charts $\{(U_{\al},\phi_{\al})\}$ for $(X,p)$
as in Definition~\ref{def:cf} such that:
\bi
\item[(i)] for each $\al$, there is given a trivialization
$(\pi^{-1}(U_{\al}),\Psi_{\al})$ of $E\mid U_{\al}$;
\item[(ii)] with $\Psi_{\al}$ as in (i) and with
$\kappa_{\al}=(\phi_{\al}\otimes 1)\Psi_{\al}:\pi^{-1}(U_{\al})\to Y\x\R^{k+m}$, every
each of the  maps $\kappa_{\beta}\kappa_{\al}^{-1}$ belongs to
$\text{Diff}^{0}_{\ell}(\kappa_{\al}(\pi^{-1}(U_{\al})),
\kappa_{\bt}(\pi^{-1}(U_{\bt})))$.
\ei
\end{definition}

It is left to the reader to check that the standard operations on vector
bundles (such as those of forming alternating and tensor products of bundles)
preserve the smoothness property.

An important smooth vector bundle over $X$ is the tangent bundle $(TX,\pi)$,
where
\[ TX = \cup_{y\in Y}TX^{y}        \]
and $\pi$ is the canonical projection map.
Let us show that $TX$ is indeed a smooth $k$-dimensional
vector bundle over $X$.  Let
$\{(U_{\al},\phi_{\al}): \al\in A\}$ be an atlas for $(X,p)$
and let
$Z_{\al}=\pi^{-1}(U_{\al})$ for each $\al$.  For each $y\in p(U_{\al})$,
the restriction map $\phi_{\al}^{y}$ of $\phi_{\al}$ to
$X^{y}_{\al}=X^{y}\cap U_{\al}$ is a diffeomorphism onto an open subset of
$\R^{k}$. Define
$\Psi_{\al}:Z_{\al}\to U_{\al}\x \R^{k}$ by:
\[   \Psi_{\al}(\ga)=(\phi_{\al}(\pi(\ga)),
D_{\pi(\ga)}\phi_{\al}^{p\circ\pi(\ga)}(\ga)).               \]
It is easily checked that each $\Psi_{\al}$ is a bijection onto a set of the
form $A\x \R^{k}$ where $A$ is an open subset of $Y\x \R^{k}$, and that every
$\Psi_{\al}(Z_{\al}\cap Z_{\bt})$ is also an open subset of this form.
The transition function
$\Psi_{\bt}\circ\Psi_{\al}^{-1}
:\Psi_{\al}(Z_{\al}\cap Z_{\bt})\to \Psi_{\bt}(Z_{\al}\cap Z_{\bt})$
is given by:
\[   (y,w,\xi)\to ((\phi_{\bt}^{y}\circ (\phi_{\al}^{y})^{-1})(y,w),
D_{w}(\phi_{\bt}^{y}\circ (\phi_{\al}^{y})^{-1})(\xi)).            \]
It is routine to check that
$\Psi_{\bt}\circ\Psi_{\al}^{-1}\in \Ga(Y\x \R^{k};\R^{k})$.  So
$TX$ is a smooth vector bundle over $X$.

Let $E$ be a smooth vector bundle over $X$.  A {\em Riemannian metric} for
$E$ is a family of Riemannian metrics on the smooth
vector bundles $E\mid X^{y}$
$(y\in Y)$ which, in terms of local coordinates, vary continuously. Using a
$C^{0,\infty}-$ partition of unity, a standard argument shows that $E$ admits a
Riemannian metric.

\section{Continuous family groupoids}
We now discuss the class of locally compact groupoids with which this paper
is primarily concerned and which generalize Lie groupoids.  We first recall
some facts about locally compact groupoids.

A {\em groupoid} is most simply defined as a small category with inverses.
Spelled out axiomatically, a groupoid
is a set $G$ together with a subset $G^{2} \subset G\times G$, a product map
$m:G^{2}\to G$, where we write $m(a,b)=ab$,
and an inverse map $i:G\to G$, where we write $i(a)=a^{-1}$
and  where $(a^{-1})^{-1}=a$,  such that:
\bi
\item[(i)] if $(a,b), (b,c)\in G^{2}$, then $(ab,c), (a,bc)\in G^{2}$
and
\[     (ab)c=a(bc);        \]
\item[(ii)] $(b,b^{-1})\in G^{2}$ for all $b\in G$, and if $(a,b)$
belongs to $G^{2}$, then
\[     a^{-1}(ab) = b \hspace{.2in} (ab)b^{-1} = a. \]
\ei
We define the {\em range} and {\em source} maps $r:G\to G^{0}$, $s:G\to G^{0}$
by:
\bgc
\[       r(x)=xx^{-1}\hspace{.2in} s(x)=x^{-1}x.     \]
\edc
The {\em unit space} $G^{0}$ is defined to be $r(G)=s(G)$, or equivalently,
the set of idempotents $u$ in $G$. The maps $r,s$ fiber the groupoid $G$ over
with fibers $\{G^{u}\}, \{G_{u}\}$, where $G^{u}=r^{-1}(\{u\})$ and
$G_{u}=s^{-1}(\{u\})$.  Note that $(x,y)\in G^{2}$ if and only if
$s(x)=r(y)$.

For detailed discussions of groupoids (including locally compact and Lie
groupoids below), the reader is referred to the books
\cite{Mackenzie,Paterson,rg}. Important examples of groupoids are given by
transformation group groupoids and equivalence relations.

A {\em locally compact groupoid} is a groupoid $G$ which is also a second
countable locally compact Hausdorff space for which multiplication and
inversion are continuous. (A detailed discussion of non-Hausdorff locally
compact groupoids is given in \cite{Paterson}.) Note that $G^{2},G^{0}$ are
closed subsets of $G\x G, G$ respectively.  Further, since $r,s$ are
continuous, every $G^{u}, G_{u}$ is a closed subset of $G$.

A locally compact groupoid $G$ is called a {\em Lie groupoid} if $G$ is a
manifold such that:
\bi
\item[(i)] $G^{0}$ is a submanifold of $G$;
\item[(ii)] the maps $r,s:G\to G^{0}$ are submersions;
\item[(iii)] the product and inversion maps for $G$ are smooth.
\ei
Note that $G^{2}$ is naturally a submanifold of $G\x G$ and every $G^{u},
G_{u}$ is a submanifold of $G$.  (See \cite[pp.55-56]{Paterson}.)

For analysis on a locally compact groupoid $G$, it is essential to have
available a {\em left Haar system}.  This is the groupoid version of a {\em
left Haar measure}, though unlike left Haar measure on a locally compact
group, such a system may not exist and if it does, it will not usually be
unique.  However, in many case, there is a natural choice of left Haar system.
For Lie groupoids, such a system exists and is essentially unique.  As we
will see later, this result extends to continuous family groupoids defined
below.

A left Haar system on a locally compact groupoid $G$ is a family of measures
$\{\la^{u}\}$ $(u\in G^{0})$, where each $\la^{u}$ is a
positive regular Borel measure on the locally compact Hausdorff space
$G^{u}$, such that the following three axioms are satisfied:
\bi
\item[(i)] the support of each $\la^{u}$ is the whole of $G^{u}$;
\item[(ii)] for any $g\in C_{c}(G)$, the function $g^{0}$, where
\[   g^{0}(u)=\int_{G^{u}}g\,d\la^{u},   \]
belongs to $C_{c}(G^{0})$;
\item[(iii)] for any $x\in G$ and $f\in C_{c}(G)$,
\[ \int_{G^{d(x)}}f(xz)\,d\la^{d(x)}(z) = \int_{G^{r(x)}}f(y)\,d\la^{r(x)}(y).\]
\ei
The existence of a left Haar system on $G$ has topological consequences for
$G$ -- it entails that both $r,s:G\to G^{0}$ are open maps
(\cite[p.36]{Paterson}).

\begin{definition}           \label{def:sfg}
A locally compact groupoid $G$ is called a {\em continuous family groupoid} if:
\bi
\item[(i)] both $(G,s), (G,r)$ are continuous families of manifolds over
$G^{0}$;
\item[(ii)] the inversion map $i:(G,s)\to (G,r)$, where $i(x)=x^{-1}$,
is an isomorphism of continuous families of manifolds:
\begin{equation}  \label{CD:inver}
\begin{CD}
	G  @>i>>         G  \\
      @VsVV         @VVrV \\
	G^{0}   @>=>>    G^{0}
\end{CD}
\end{equation}
\item[(iii)] the multiplication map $m:(G*G,t_{1})\to (G,r)$ is a morphism of
continuous families with respect to $r$:
\begin{equation}  \label{CD:G*G}
\begin{CD}
	G*G  @>m>>         G  \\
      @Vt_{1}VV              @VVrV \\
	G       @>r>>       G^{0}
\end{CD}
\end{equation}
\ei
\end{definition}

In (iii) above, $G*G$ is defined to be $(G,s)*(G,r)$, so that
\[      G*G=\{(x,y)\in G\x G: s(x)=r(y)\} = G^{2}.     \]
As in \S 2, $G*G$ is a continuous family over $G$ with $t_{1}((x,y))=x$.
The map $m$ is a fiber preserving map from $(G*G,t_{1})$ into $(G,r)$ (with
respect to $r:G\to G^{0}$) since $r(xy)=r(x)$.

Condition (iii) as stated
is one-sided in that the morphism property of $m$ is
formulated in terms of $t_{1}, r$ rather than  $t_{2},s$.  Indeed, there is a
commuting diagram:
\begin{equation}  \label{CD:G*Gs}
\begin{CD}
	G*G  @>m>>         G  \\
      @Vt_{2}VV              @VVsV \\
	G       @>s>>       G^{0}
\end{CD}
\end{equation}
and it is as natural to formulate the morphism property for $m$ in terms of
(\ref{CD:G*Gs}) as in terms of (\ref{CD:G*G}). We now show that these two
formulations are actually equivalent given (i) and (ii) of
Definition~\ref{def:sfg}.

\begin{proposition}        \label{prop:symm}
Let $G$ be a locally compact groupoid satisfying (i) and (ii) of
Definition~\ref{def:sfg}.  Then the multiplication map
$m:(G*G,t_{1})\to (G,r)$ is a morphism with
respect to $r$ if and only if $m:(G*G,t_{2})\to (G,s)$ is a morphism with
respect to $s$.
\end{proposition}
\begin{proof}
We show first that $(i*i)^{\Tilde{}}$, where
$(i*i)^{\Tilde{}}((x,y))=(y^{-1},x^{-1})$, is a morphism from $(G*G,t_{1})$
to $(G*G,t_{2})$ with respect to $i$:
\begin{equation}        \label{CD:i*i}
\begin{CD}
     G*G  @>(i*i)^{\Tilde{}}>>     G*G \\
	@Vt_{1} VV           @Vt_{2}VV \\
	G           @>i>>         G
\end{CD}
\end{equation}
Clearly, $(i*i)^{\Tilde{}}$ is fiber preserving with respect to $i$.
Let $(x,y)\in G*G$.  Let $U'$ be a chart for $(G*G,t_{2})$ containing
$(y^{-1},x^{-1})$. We can suppose that $U'=U_{1}*U_{2}$ where $U_{1}$ is a
chart for $(G,s)$, $U_{2}$ is a chart for $(G,r)$ and $s(U_{1})=r(U_{2})$.
We can suppose further that
$U_{1}\sim s(U_{1})\x W_{1}, U_{2}\sim r(U_{2})\x W_{2}$. Then
$U_{1}*U_{2}\sim U_{2}\x W_{1}$ under the map
$(u_{1},u_{2})\to (u_{2},w_{1})$
where $u_{1}\sim (s(u_{1}),w_{1})$. Note that
$((i*i)^{\Tilde{}})^{-1}(U_{1}*U_{2})=U_{2}^{-1}*U_{1}^{-1}$.

By (ii) of Definition~\ref{def:sfg}, there exists a chart $U_{1}'$ for $(G,s)$
with $x\in U_{1}'$ and $U_{1}'\sim s(U_{1}')\x W_{1}', i(U_{1}')\subset
U_{2}$. Similarly, there exists a chart $U_{2}'$ for $(G,r)$
with $y\in U_{2}'$ and $U_{2}'\sim r(U_{2}')\x W_{2}', i(U_{2}')\subset
U_{1}$. Since $s(x)=r(y)$, we can, by contracting $U_{1}', U_{2}'$,
suppose that $s(U_{1}')=r(U_{2}')$. Then $(x,y)\in U_{1}'*U_{2}'\sim U_{1}'\x
W_{2}'$ is a chart for $(G*G,t_{1})$.  Now in local terms,
for $(x_{1},w_{2}')\in U_{1}'\x W_{2}'$,
\[   (i*i)^{\Tilde{}}(x_{1},w_{2}')=
(i^{-1}(s(x_{1}),w_{2}'),x_{1}^{-1})  \]
so that $(i*i)^{\Tilde{}}$ is a morphism, again using (ii) of
Definition~\ref{def:sfg}.  

Suppose that $m:(G*G,t_{2})\to (G,s)$ is a morphism with respect to $s$.
Then $m:(G*G,t_{1})\to (G,r)$ is a morphism since it is the morphism product
$i\circ m\circ (i*i)^{\Tilde{}}$:
\[
\begin{CD}
     G*G  @>(i*i)^{\Tilde{}}>> G*G   @>m>>     G       @>i>>         G \\
     @Vt_{1} VV             @Vt_{2} VV         @VsVV                 @VVrV  \\
      G          @>i>>        G      @>s>>   G^{0}  @>=>>   G^{0}
\end{CD}
\]
The converse is true in a similar way.
\end{proof}

Let $G$ be a continuous family groupoid. Note that $(G*G)^{x}=\{x\}\x G^{s(x)}$
in $(G*G,t_{1})$, so that what (iii) of Definition~\ref{def:sfg} above is
saying is that for fixed $x\in G$, the map $L_{x}:y\to xy$ is a diffeomorphism
from $G^{s(x)}$ onto $G^{r(x)}$ and this diffeomorphism is required to vary
continuously with $x$. The same holds for the right multiplication map
$R_{y}:x\to xy$ by Proposition~\ref{prop:symm}. So the multiplication in $G$
can be regarded as ``separately continuous'' in the sense that the
diffeomorphisms $L_{x},R_{y}$'s vary continuously. In contrast, Lie groupoids
satisfy the much stronger ``joint continuity'' condition that the
multiplication map $m:G*G\to G$ is smooth.

It is convenient to have available the pull-back version of (\ref{CD:G*G}).
Note that the pull-back $r^{-1}G$ of $G,r)$ is:
\[   \{(x,z)\in G\x G: r(x)=r(z)\}. \]
The latter set will be denoted by $G*_{r} G$ and is the fibered product
$(G,r)*(G,r)$. (Similarly, $G*_{s}G$ is the fibered product
$(G,s)*(G,s)$.) The associated morphism $m':G*G\to
G*_{r} G$ is given by: $m'(x,y)=(x,xy)$ and we have the diagram:
\begin{equation} \label{CD:morgpd}
\begin{CD}
	G*G      @>m'>>  G*_{r} G   \\
  @Vt_{1}VV             @VVt_{1}V \\
	 G       @>=>>      G
\end{CD}
\end{equation}

We now discuss some examples of continuous family groupoids. Firstly, every Lie
groupoid is a continuous family groupoid. Indeed, in Definition~\ref{def:sfg},
(i) follows since $r,s$ are submersions. Properties (ii) and (iii) follow since
both $i$ and $m$ are smooth.

A very simple example of a continuous family groupoid that is not a Lie
groupoid is provided by any locally compact Hausdorff space $Y$ that is not a
manifold (treated as a groupoid of units). This is the groupoid associated
with the (non-equivariant) index theorem for families.

Next, the transformation group groupoid $G$ given by an action of a Lie group
$H$ on a locally compact Hausdorff space $Y$ is a continuous family groupoid.
In this case, $G=H\x Y$, and the multiplication is given by
$((h,y),(k,k^{-1}y))\to (hk,k^{-1}y)$ and inversion by $(h,y)\to (h^{-1},hy)$.
The unit space $G^{0}$ can be identified with $Y$ and the source and range maps
are given by: $s((h,y))=y, r((h,y))=hy$. For $y\in Y$, we have $G_{y}=H$,
$G^{y}= \{(h,h^{-1}y):h\in H\}$ which can also be identified with $H$ by
sending $(h,h^{-1}y)$ to $h$. We will not give the proof here that $G$ is a
continuous family groupoid since this will be generalized later in
Proposition~\ref{prop:gy}. It is easy to see that $G$ is a Lie groupoid if and
only if $Y$ is a manifold on which $H$ acts smoothly. This gives many examples
of continuous family groupoids which are not Lie groupoids. For example, if $H$
acts trivially on $Y$ and $Y$ is not a manifold, then $G$ is a continuous
family groupoid that is not a Lie groupoid.

Another example of a continuous family groupoid is an equivalence relation
$G=\{((y,z),(y,z')): y\in Y, z\in Z\}$ on $Y\x Z$, where $Z$ is a smooth
manifold. So $G=Y\x (Z\x Z)$ with the product given by
$(y,z,z')(y,z',z'')=(y,z,z'')$ and inversion by $(y,z,z')\to (y,z',z)$. The
unit space of $G$ is $Y\x Z$, where $Z$ is identified with the diagonal
$\{(z,z):z\in Z\}$ in the obvious way. The source and range maps are given by:
$s(y,z,z')=(y,z'), r(y,z,z')=(y,z)$. Each of $G^{(y,z)}, G_{(y,z)}$ is just
the manifold $Z$. It is easy to check that $G$ is a continuous family groupoid,
and is a Lie groupoid if and only if $Y$ is a manifold.

Let $G$ be a continuous family groupoid.  From the discussion in \S 2 (with
$(X,p)=(G,r)$), $TG=\cup_{u\in G^{0}}TG^{u}$
is a smooth vector bundle over $(G,r)$. The  restriction $A(G)$ of $TG$ to
$G^{0}$ is a vector bundle over $G^{0}$, and, as in the case of Lie
groupoids, is called the {\em Lie algebroid} of $G$.  We now briefly discuss
the existence of left Haar systems on $G$. The discussion parallels that for
Lie groupoids given in \cite[2.3]{Paterson}.

For a chart $(U,\phi)$ of (the continuous family) $(G,r)$ and for any
$u\in r(U)$, let $\phi^{u}$ be the restriction of $\phi$ to $U\cap G^{u}$
and $W^{u}=\phi^{u}(U\cap G^{u})\subset \R^{k}$ (where $k$ is the dimension
of the manifolds $G^{u}$).
Given a measure $\mu$ on $U\cap G^{u}$,
the  measure $\mu\circ \phi^{u}$ on $W^{u}$ is defined by:
\[         \mu\circ \phi^{u}(E)=\mu((\phi^{u})^{-1}(E)).         \]
For any open subset $W$ of $\R^{k}$
let $\la^{W}$ be Lebesgue measure restricted to $W$.
A left Haar system $\{\la^{u}\}$ for $G$ is called a
{\em $C^{0,\infty}$ left Haar
system} if for any chart $(U,\phi)$ for $(G,r)$, the measure
$\la^{u}\circ \phi^{u}$ is equivalent to $\la^{W^{u}}$, and the map $f$ that
sends
$(u,t)\to (d(\la^{u}\circ \phi^{u})/d\la^{W^{u}})(u,t)$ belongs to
$C^{0,\infty}(\phi(U))$.  The left  Haar system is called {\em continuous} if
$f\in C(\phi(U))$.   (Of course every $C^{0,\infty}$ left Haar
system is continuous.)  If $G$ is a Lie groupoid, then the left Haar system
is called {\em smooth} if $f$ is a $C^{\infty}$ function.

Note that two continuous left Haar systems on a continuous family groupoid give
isomorphic universal $\css$ and isomorphic reduced $\css$. The argument for
this (which presupposes Renault's representation theory of locally compact
groupoids (\cite{MuhlyTCU,Ren87}) goes as follows. Any quasi-invariant
measure $\mu$ on $G^{0}$ determines, for any left Haar system $\{\la^{u}\}$, a
measure $\nu=\int\la^{u}\,d\mu(u)$ on $G$, and this in turn determines other
measures $\nu^{-1}, \nu^{2}$ on $G$ and $G^{2}$ respectively. Two continuous
left Haar systems will give equivalent measures $\nu$ (and similarly for
$\nu^{-1}, \nu^{2}$) since both $\nu$'s are equivalent on a chart. The
representations for $G$ are then the same for each of these two left Haar
systems and any such representation gives equivalent representations for each
system on $C_{c}(G)$. It follows that the universal \css\ $C^{*}(G)$ for the
two systems are isomorphic. Similarly, the reduced \css\ $C^{*}_{red}(G)$ are
also independent of the choice of system.

In the Lie groupoid case, the equivalence of smooth left Haar systems is
captured precisely in Connes's density bundle approach to integration on Lie
groupoids (\cite[p.101]{Connesbook}) which is canonical, and indeed that
approach can be readily adapted to apply to $C^{0,\infty}$
(resp. continuous) left
Haar systems for continuous family groupoids. However, for relating the
representation theory for continuous family groupoids to the existing
representation theory for locally compact groupoids as well as for calculation
purposes, the $C^{0,\infty}$ (resp. continuous) left Haar system approach is
convenient and will be used in this paper.

It is known  that every Lie groupoid admits a smooth
left Haar system.  The present writer does not know if every continuous
family groupoid admits a $C^{0,\infty}$ left Haar system. However, as we
shall see, every continuous family groupoid admits a {\em continuous} left Haar
system. The proof of this result is along the same lines as that for the
existence of a smooth left Haar system on a Lie groupoid
(\cite[p.63]{Paterson}).  (See also \cite{Landsman,Ramazan}.)

\begin{theorem}  \hspace{.1in}        \label{th:LHS}
Let $G$ be a continuous family groupoid.  Then there exists a continuous left
Haar system on $G$.
\end{theorem}
\begin{proof}
We observe first that the $1$-density bundle $\Om(A(G)^{*})$
is trivial, and there exists a strictly positive section $\al$ of that bundle.
For each $x\in G$, define $L_{x^{-1}}:G^{r(x)}\to G^{s(x)}$ by:
$L_{x^{-1}}(y)=x^{-1}y$. The map
$L_{x^{-1}}$ is a diffeomorphism (using Definition~\ref{def:sfg}).

As in the Lie groupoid case, we take the measure $\la^{u}$ to be the regular
Borel measure associated with the density $z\to
(L_{z^{-1}})^{*}_{s(z)}(\al_{s(z)})$ on $G^{u}$.  To check that this is well
defined, for each $u_{0}\in G^{0}$ and $z_{0}\in G^{u_{0}}$, we
obtain in terms of local charts for $u_{0}, z_{0}$ in $G$ that
$\al_{u}=g(u)\,dw_{1}\ldots dw_{k}$ for some continuous positive function $g$,
and
\[  (L_{z^{-1}})^{*}_{s(z)}(g(s(z))\,dw_{1}\cdots dw_{k})
= g(s(z))\mid J(L_{z^{-1}})(z)\mid dz_{1}\cdots dz_{k}  \]
where $J$ stands for the Jacobian.
The function  $z\to g(s(z))\mid J(L_{z^{-1}})(z)\mid$ is (using
Definition~\ref{def:sfg}) a continuous function.
It follows that $\{\la^{u}\}$ is a continuous left Haar system
on $G$, the proofs of the other items requiring to be checked
being the same as in the Lie groupoid case.
\end{proof}

\section{Actions of continuous family groupoids}
To motivate the need for continuous family groupoid actions, it is helpful to
consider the situation of
the Atiyah-Singer equivariant families theorem.  (Atiyah and Singer refer
this ``to the reader'' (\cite[p.135]{AS4}).) There, we have a compact
fiber bundle $(X,p)$ over $Y$ with compact smooth manifold $Z$ as fiber and
structure group $\text{Diff}(Z)$.  We are also given a compact Lie group $H$
acting continuously on $Y$ and in a $C^{0,\infty}$-way on $X$ with $p$
equivariant.  We want to interpret this in terms of an action of the
continuous family groupoid $H*Y$ on $X$, given by:
\[            (h,y)x=hx      \]
where $p(x)=y$.  The index of an equivariant family of pseudodifferential
operators elliptic along the leaves will then lie in $K_{0}(C^{*}(H*Y))$.
(The Lie groupoid version of this is proved in \cite{Paterson2}.)  So
even in the classical Atiyah-Singer index context, we have to leave the
category of groups and move to the category of continuous family groupoids.

With this motivation, we now turn to the problem of determining what a left
action of a continuous family groupoid should be. It turns out - slightly
surprisingly - that all that
we require for present purposes is a {\em continuous} action
of a continuous family groupoid $G$ on a fiber space $Y$ over $G^{0}$. The
reason is that the unit space of the locally compact groupoid $G*Y$ is
identified with $Y$, while its fibers $(G*Y)_{y}, (G*Y)^{y}$ are effectively
the continuous family groupoid
fibers $G_{p(y)}, G^{p(y)}$, so that the $C^{0,\infty}-$
structure of $G*Y$ is effectively that of $G$ with the $Y$ only playing the
role of a continuous parameter space.  We now discuss this in more detail.

Let $G$ be a locally compact groupoid (not necessarily a continuous family
groupoid) with $r,s$ open maps, and let
$(Y,p)$ be a fiber space over $G^{0}$. Form
the fibered product $G*Y=(G,s)*(Y,p)$ of the fiber spaces $(G,s), (Y,p)$
over $G^{0}$. The action of $G$ on $Y$ is then given by a continuous map
$n:G*Y\to Y$. The action has to satisfy the natural algebraic axioms: so we
require $p(gy)=r(g)$, $g_{1}(g_{2}y)=(g_{1}g_{2})y$ and $g^{-1}(gy)=y$ whenever
these make sense. The space $Y$ with such an action of $G$ is called a
{\em $G$-space} (cf.\cite{MuReWi}).

Let $Y$ be a $G$-space. It is well known and easy to check from the groupoid
axioms of \S 2 that $G*Y$ is a locally compact groupoid with
operations given by: $(h,z)(g,g^{-1}z)=(hg,g^{-1}z)$ and
$(g,y)^{-1}=(g^{-1},gy)$. Let $\mu$ be the multiplication map and $i$ be
the inversion map on $G*Y$. So
\begin{equation}  \label{eq:mu}
  \mu(h,z,g,g^{-1}z)=(hg,g^{-1}z).
\end{equation}
Further, $s(g,y)=(s(g),y), r(g,y)=(r(g),gy)$. We can identify
$(s(g),y)$ with $y$ (since $s(g)=p(y)$), and so
$(G*Y)^{0}=Y$. (In fact, this is just identifying $G^{0}*Y$ with $Y$.)
Note that $r=n$ and $s=t_{2}$.

Now, in addition, assume that $G$ is a continuous family groupoid. We
will show that $G*Y$ is also, in a natural way, a continuous family groupoid.
We first discuss the continuous family structures for $(G*Y,s)$ and
$(G*Y,r)$. We have a commuting diagram:
\begin{equation}  \label{CD:beta}
\begin{CD}
	G*Y      @>i>>      G*Y  \\
      @VrVV                @VVsV \\
	Y       @>=>>        Y
\end{CD}
\end{equation}
where $i$ is the inversion map on $G*Y$.
We can
identify $(G*Y,s)$ canonically with $(G*Y,t_{2})$.  As observed in \S 2,
$(G*Y,t_{2})$ is a continuous family (since $(G,s)$ is a continuous family).
Next, the map $r$ is a continuous surjection which is open,
since $s$ is and $i$
is a homeomorphism. We give $(G*Y,r)$ the (unique) continuous family structure
that makes $i$ an isomorphism of continuous families. In particular, an atlas
for $(G*Y,r)$ is determined by charts of the form $i^{-1}(U*V)$, and the
$C^{\infty}$-structure on each $(G*Y)^{z}=\{(g,g^{-1}z): g\in G^{p(z)}\}$ is
just that obtained by identifying $(G*Y)^{z}$ with $G_{p(z)}$ through the map
$(g,g^{-1}z)\to g^{-1}$ (and hence with $G^{p(z)}$ through the map
$(g,g^{-1}z)\to g$, using (ii) of Definition~\ref{def:sfg}.)

\begin{proposition}  \hspace{.1in}        \label{prop:gy}
The continuous family $G*Y=(G,s)*(Y,p)$ is a continuous family groupoid.
\end{proposition}
\begin{proof}
We check that the requirements for a continuous family groupoid in
Definition~\ref{def:sfg} hold for $G*Y$. (i) and (ii) of that definition
follow from the discussion preceding the statement of the proposition.
It remains to show that the multiplication map $\mu$ on $G*Y$ is a morphism.
To this end, consider the following commutative diagram:
\begin{equation}  \label{CD:tilmu}
\begin{CD}
	(G*Y)*(G*Y)      @>\al>> G*Y  @>i^{-1}>>     G*Y  \\
      @Vt_{1}VV         @VsVV                    @VrVV \\
	  G*Y    @>r>>    Y    @>=>>        Y
\end{CD}
\end{equation}
Here, $\al(h,z,g,g^{-1}z)=(g^{-1}h^{-1},hz)$.  Then $\mu=i^{-1}\circ\al$,
and since $i^{-1}$ is a morphism ((\ref{CD:beta})), it remains to show that
$\al$ is a morphism.  It is sufficient for this to show that $\ga:G*G\to G$,
where $\ga(h,g)=g^{-1}h^{-1}$, is a morphism:
\begin{equation}
\begin{CD}  \label{CD:GGG}
		G*G     @>\ga>>     G  \\
	     @Vt_{1}VV           @VVsV  \\
		 G      @>r>>       G^{0}
\end{CD}
\end{equation}
For since the map $(h,z)\to hz$ is continuous, it follows that
$(h,z)\to \al_{(h,z)}$ varies in a $C^{0,\infty}-$ way if the map
$h\to \ga_{h}$ does.  The map $\ga$ is a morphism since it is the composition
of two morphisms:
\begin{equation}  \label{CD:GGGG}
\begin{CD}
	G*G      @>m>>    G       @>i>>     G  \\
      @Vt_{1}VV         @VrVV             @VsVV \\
	  G      @>r>>  G^{0}     @>=>>    G^{0}
\end{CD}
\end{equation}
\end{proof}

Motivated by the families situation  (discussed in the
first paragraph of this section) we now have to extend the $G$-space notion to
that of a fiber space $(X,q)$ over $Y$ which is also a $G$-space with an
action compatible with that on $Y$.
compatible action.  Precisely, $(X,p\circ q)$ is a fiber space over $G^{0}$
and we require that this fiber space and $(Y,p)$ be  $G$-spaces such that
for all $(g,x)\in G*X$, we have $q(gx)=gq(x)$.  (Note that $(g,q(x))\in G*Y$
since $p(q(x))=s(g)$.)  This can be formulated in terms of the commuting
diagram of continuous  maps:
\begin{equation}  \label{CD:GYX}
\begin{CD}
	G*X  @>n_{X}>>    X  \\
      @V\text{id}*qVV         @VVqV \\
	G*Y   @>n_{Y}>>   Y
\end{CD}
\end{equation}
where $n_{X}, n_{Y}$ are the action maps of $G$ on $X$ and $Y$ respectively.
We say that {\em $X$ is a $G$-space over $Y$}.
We now show that $X$ is itself a $G*Y$-space in a natural way.

\begin{theorem}   \label{th:final}
Let $G$ be a continuous family groupoid and $Y$ be a $G$-space.  Then the
class of $G$-spaces $X$ over $Y$ is canonically identified with the class of
$G*Y$-spaces.
\end{theorem}
\begin{proof}
Let $(X,q)$ be a $G$-space over $Y$. Recalling that the source map
of the groupoid
$G*Y$ is the map $(g,y)\to y$, it follows that $(G*Y)*X=(G*Y,s)*(X,q)$ is a
fiber space over $Y$. Define a map
\begin{equation}  \label{eq:nYX}
    n_{Y,X}:(G*Y)*X\to X,
\end{equation}
by: $n_{Y,X}(g,y,x)=gx$. Calculations very similar to those of the next
paragraph show that $n_{Y,X}$ is an action of $G*Y$ on $X$. So $X$ is a
$G*Y$-space.

Conversely, suppose that $(X,q)$ is a $G*Y$-space.
Then $q:X\to (G*Y)^{0}=Y$ is continuous and onto,
and so $X$ is a fiber space over $Y$. Let
$n:(G*Y)*X\to X$ be the action map.  Then $(X,p\circ q)$ is a fiber space
over $Y$, so that $G*X$ is defined.
For $(g,x)\in G*X$, define the action of $G$ on $X$ by:
       \[          gx=n((g,q(x)),x)=(g,q(x))x.   \]
Clearly the map $(g,x)\to gx$ is continuous. We now check the algebraic action
axioms for the map $(g,x)\to gx$. For the associative law, with $z=q(x)$, we
have $(hg)x=[(h,gz)(g,z)]x=(h,gz)[(g,z)x]=(h,gz)(gx)=h(gx)$. Next
$g^{-1}(gx)=s(g)x=s((g,z))x=x$. Lastly, $q(gx)=q((g,z)x)=r((g,z))=gz$. It
follows that $X$ is a $G$-space over $Y$. It is left to the reader to show that
if we apply the $n_{Y,X}$ construction to $X$ with this $G$-space structure,
then we get back to the $G*Y$-space with which we started.
\end{proof}

In the situation of the equivariant families index theorem, one has to consider
a $G$-space $X$ over $Y$ which is a continuous family and on which $G$ acts
in a $C^{0,\infty}-$ way.  We now briefly indicate how this is defined.  By
Theorem~\ref{th:final}, we can, by replacing $G$ by $G*Y$, suppose that $X=Y$.
Then the action of $G$ on $Y$ is said to be {\em $C^{0,\infty}$} if the
multiplication map $n:(G*Y,t_{1})\to Y$ is a morphism with respect to $r$:
\begin{equation}  \label{CD:G*Y}
\begin{CD}
	G*Y  @>n>>         Y  \\
      @Vt_{1}VV              @VVpV \\
	G       @>r>>       G^{0}
\end{CD}
\end{equation}
We say that $Y$ is a {\em $C^{0,\infty}$ $G$-space}.
The continuous family groupoid $G$ is itself a $C^{0,\infty}$ $G$-space.
This follows from  (iii) of Definition~\ref{def:sfg}.

As in (\ref{CD:morgpd}), the morphism property of (\ref{CD:G*Y})
can be reformulated in terms of a morphism $n'$:
\begin{equation} \label{CD:moract}
\begin{CD}
	G*Y     @>n'>>   G*_{r} Y   \\
      @Vt_{1}VV @VVt_{1}V \\
	     G @>=>> G
\end{CD}
\end{equation}
where $G*_{r} Y=(G,r)*(Y,p)$ and the map $n'$ is given by:
$n'(g,y)=(g,gy)$.

In conclusion, Theorem~\ref{th:final} says that a $G$-space $X$ over $Y$ is the
same as a $(G*Y)$-space, the fibering of $X$ over $Y$ being ``absorbed'' as it
were into the fibering of the continuous family groupoid $G*Y$. (The same
applies if we work in the category of $C^{0,\infty}$ $G$-spaces.) Note
that this cannot be formulated if we stay within the group category, since in
forming $G*Y$, we leave the group category. So we do not need
to work in the situation where a $G$-space $X$ is fibered over another
$G$-space $Y$. For the ``higher order'' fibered space $X$ is itself just an
``ordinary'' groupoid space for the groupoid $G*Y$, and we do not leave the
category of continuous family groupoids by forming $G*Y$. Thus we only ever
need consider the action of a continuous family groupoid on a $G$-space,
changing the groupoid if necessary.


\end{document}